\documentclass[12pt]{article}

\usepackage[utf8]{inputenc}
\usepackage{amssymb}
\usepackage{amsmath}
\usepackage{amsthm}
\usepackage[cjk]{kotex}
\usepackage[colorlinks]{hyperref}
\usepackage{multicol}

\usepackage{fullpage}

\theoremstyle{theorem}

\theoremstyle{definition}

\title{Doing Math in Jest:\ Reflections on Useless Math, the Unreasonable Effectiveness of Mathematics, and the Ethical Obligations of Mathematicians}
\author{Gizem Karaali}

\begin{document}

\date{}

\maketitle

\begin{abstract}{Mathematicians occasionally discover interesting truths even when they are playing with mathematical ideas with no thoughts about possible consequences of their actions. This paper describes two specific instances of this phenomenon. The discussion touches upon the theme of the unreasonable effectiveness of mathematics as well as the ethical obligations of mathematicians.}\end{abstract}

\section*{Mathematicians at Play:\ Doing Math in Jest}

In 1993, Jean-Fran\c{c}ois Mestre, Ren\'{e} Schoof, Lawrence Washington, and Don Zagier published ``Quotients homophones des groupes libres / Homophonic Quotients of Free Groups'' in the then-new journal {\em Experimental Mathematics} (volume {\bf 2} number {3}, pages 153--155) \cite{QHGL}. This was a joke of a paper; the authors were simply playing around. In one half of the paper, written in French, the authors took the English alphabet, defined the free group generated by it, and then looked for relations defined by the homophones in the language to quotient this behemoth with. Taking the reader through identities such as ``lead = led", ``knight = night", and ``dammed = damned", they concluded that the quotient is trivial. The second half of the paper, written in English, did the same thing to French. So in this four-author, three-page, two-column paper, the main result was that both English and French were trivial; or rather, the homophonic quotient groups of English and French were trivial. The fun continued into the acknowledgments, where the third and fourth authors mentioned being partially supported by the NSF, ``and (by the results of this paper) numerous other government agencies". (Here is a start: NSF = NSA = NASA = NAS = NAE = NEA = NEH = NIH.)

Then in 2016, Herbert Gangl, at Durham University, Woohyung Lee, now at University of Illinois Chicago, and I got together to discuss and explore what would happen to other languages we were familiar with. Herbert had seen the 1993 paper earlier, and had already mostly worked out the German analogue of the problem; he thought the problem was fun. For Woohyung and myself, this puzzle was also a neat distraction. We were reading one of the yellow Springer GTM volumes together, and whenever we got burnt out, we switched to playing around with words in Korean (his native language) and Turkish (mine). Eventually, Woohyung reached out to Herbert, and we decided to put our results together and submit the whole project for publication. 

The editors of the first journal we submitted our paper to told us that this was not ``serious research". It was also not as ``entertaining" as the original 1993 paper had been, so it was not acceptable for this journal. But don't many papers have histories of prior rejection? This one, too, finally found a home for itself \cite{GKL}. 

But I am not writing to bemoan the idiosyncracies of academic publishing. This is about what our paper and its initial rejection helped me see:\ even when math is done in jest, sometimes interesting truths come out of it. The authors of \cite{QHGL} were having fun; Herbert, Woohyung and I, too, had a lot of fun working on the project. But we were also able to uncover a nugget of truth through this mathematical play. 

\section*{The Unreasonable Effectiveness of Mathematics}

Readers probably recognize immediately that I am reaching toward that ineffable construct, the unreasonable effectiveness of mathematics. Introduced by Eugene Wigner in \cite{Wigner}, this phrase captures the fascinating capacity of mathematical ideas to find unexpected applications. Though Wigner's original assertions mainly related to the applications of mathematics in the physical sciences, others after him explored the power and limitations of mathematical ideas in such diverse areas as economics \cite{Velu}, ecology \cite{And, Ginzburg}, and molecular biology \cite{Lesk1, Lesk2}. 

Many philosophers and others who think much more deeply and clearly than I have written about this theme and more generally about just how mathematics comes to have applications to the world, see for instance \cite{Colyvan,
Ferreir, Grat, Nicholson, Steiner} for an eclectic selection of relevant reflections. 
And I do not presume to have a lot to add to this sophisticated conversation (cf.\ \cite{KarUR}). However, I do want to propose our little piece of mathematics done in jest as a case study, as an instance to explore. 

Our problem was not urgent, it was not intended to solve a global or even a local issue, and it was certainly not originally conceived of as applied math. We were simply curious to see what we would get when we plugged our individual languages into the basic construction proposed in \cite{QHGL}. Once we saw that German yielded a trivial group, while the Korean group was nontrivial but much smaller than the original free group, and Turkish did not have much to quotient with unless we pushed hard into ``how words are actually pronounced by real live people" \cite{Wiegold}, we thought a bit about what all this really meant. 

We finally converged to the following:\ the homophonic quotient groups offer a crude way to measure the phonetical representativeness of the associated written scripts. In other words, these groups basically measure how far the alphabet-based writing system of a specific language varies from being one-to-one representative of the sounds of the spoken version of the language. ``Languages evolve, and oral traditions evolve much faster than written ones. Thus a young script like Modern Turkish might be naturally more representative of the phonetical structure of the language and equivalently offer fewer homophones than a script which is more mature, such as the Korean one, which in turn may offer fewer homophones than an even older script such as the German one" \cite{GKL}. This is corroborated by the results of \cite{QHGL}, as both English and French are rather mature scripts.

What we had was a ``neat example of applied algebra" \cite{GKL}, adding to an already rich collection of applications of group theory to natural, sociological, and anthropological contexts. It was not a perfect measuring tool, it was not sophisticated math, and it certainly was not ``serious research". But even so, it let us verify one more time that math can help us see something about ourselves, about our human condition.

Isn't this what we humans are great at? We start with a single sentence (think ``Let there be light" or ``It was the best of times" or ``We take these truths to be self-evident") and go from there to the highest peaks and the deepest sinkholes and create volumes and volumes of scholarship, complicating and enriching the sentence ad infinitum. If we are clever and creative, we can find all sorts of ways in which our special little sentence relates to our emotional condition, the weather patterns, the last election results, and so on. Transferring the analogy to mathematics, can't we find some application for any math if we stare at it long enough?

\section*{Useless Mathematics}

G.H. Hardy famously wrote that ``a science is said to be useful if its development tends to accentuate the existing inequalities in the distribution of wealth, or more directly promotes the destruction of human life". And digging into his writing, we can see that he ``took pride in the uselessness of his work because it freed him from contributing to the terrors of war and violence" \cite{Garcia}. Even if ultimately the esoteric number theory he worked on as well as the sideways letter he threw at the young science of genetics both found themselves serious and significant applications \cite{Garcia}, many mathematicians today share this sentiment. Doing pure math feels clean because we are not doing work that directly leads to more people dying of hunger, more people being bombed, more people making money off of chemicals they keep inaccessible to those who will die without them, and so on. 

Here is another example of math done in jest, that might appear totally useless. In 2009 Philip Munz, Ioan Hudea, Joe Imad, and Robert J. Smith?\ published the first mathematical model of a zombie outbreak \cite{MHIS}. Published in a volume on infectious disease modelling, this paper uses mathematical modelling to address a fictional scenario. Even if you love zombie movies, you are probably not expecting a zombie apocalypse.\footnote{ If for some reason you are worried, the United States Centers for Disease Control and Prevention can help! See their page on Zombie Preparedness at \url{https://www.cdc.gov/cpr/zombie/index.htm}, last accessed on December 20, 2018.} The drooling zombie mobs hunting for brains are fictional. The authors approach the zombie phenomenon, all the rage in that first decade of the twenty-first century, as a context to do some interesting epidemic modelling. The zombie apocalypse mainly serves to raise the coolness factor of their work. 

I am not above using popular culture to make things fun and relevant. For a few years I taught at my institution a first-year seminar titled ``Can Zombies do Math?" \cite{K1}. In this course we used zombies to connect with the human nature of mathematics \cite{K2, K3}. 
And more recently I enjoyed reading Colin Adams' {\it Zombies and Calculus} \cite{Adams}. So in 2010 when I first came across \cite{MHIS}, I was genuinely delighted. I wanted to share the paper with everyone I knew and say: ``See! Mathematics of zombies! How cool!" 

\section*{When is Useless Math Really Useless?} 

Then I began reading. The abstract raised the first red flag: ``only quick, aggressive attacks can stave off the doomsday scenario: the collapse of society as zombies overtake us all". This sounded alarming, but I remained cheerful. 

Then I read the rest of the paper. I learned that the method proposed to save society was called ``impulsive eradication" and it required frequent attacks with increasing force, aiming to destroy more and more zombies. This seemed violent and irreversible. Other options like quarantines and cures were briefly considered, too, but these possibly more humane methods were all eventually dismissed in favor of impulsive eradication.

It felt like I was watching the authors play a first-person shooter game over their shoulders. The goal was to save the non-zombie population, so all zombies were to be killed and destroyed; the more we did this the better off we would be. First-person shooter games are fictional, too. Though people worry about their effects on children \cite{Montag}, most players can distinguish between real life and fiction. We like to think that they do know not to go on wild rampages and killing sprees in real life. 

So why worry? The paper was written in jest, its context is fictional, it is fun but irrelevant to anything else. Except possibly pedagogy. In the conclusion the authors write:\ ``[w]hile the scenarios considered are obviously not realistic, it is nevertheless instructive to develop mathematical models for an unusual outbreak. This demonstrates the flexibility of mathematical modelling and shows how modelling can respond to a wide variety of challenges in `biology'."  I can see \cite{MHIS} being used in the classroom; it could make a good reprieve from standard modelling exercises.

But this does not end here. Because the authors were quite aware of what their work implied: ``this is an unlikely scenario if taken literally, but possible real-life applications may include allegiance to political parties, or diseases with a dormant infection" \cite{MHIS}. 

\section*{Ethical Responsibilities of Doing Mathematics}

Applied mathematicians know this well:\ all models are wrong, some models are useful \cite{Box}. We knew when we submitted \cite{GKL} that our way of measuring the disparities between the spoken version of a living language and the written script that aims to represent it was far from perfect. Similarly authors of \cite{MHIS} are well aware that their model is simple and based on assumptions that make the problem easier to solve even if they might make the situation more unrealistic or contradict certain values. Still I wonder if there is something to worry about here. 

\newpage

The mathematics community sells its services to the public with propaganda slogans of the likes of ``math gives you absolute truths". Yet any mathematician working with models knows they are simplified, they are inaccurate, they are basically wrong. But we also know that there is value in them if we can get some interesting results. That is enough for the mathematician. It is the engineer, the entrepreneur, the policymaker who will need to decide how much a model applies to their context, and what to make of the related results. 

But there lies the catch. When a peer-reviewed mathematics article asserts that 
the only way to effectively combat undesirable ``allegiance(s) to political parties, or diseases with a dormant infection" is ``impulsive eradication", that quaranteens and cures won't work, what should policymakers, who may or may not be exceptionally mathematically literate, do? 

I hope it is clear that I am not advocating not doing the math. I enjoy math done in jest as much as any mathematician. Most recently I have been having all sorts of fun reading about mad veterinarians \cite{abrams}, chicken nuggets \cite{chap}, and the Game of Thrones \cite{beve}.  But I have questions. Do we have any idea how others are reading our works? And do we care at all what others might end up doing with our work? 

\section*{Useful Mathematicians Doing Useless Mathematics}

Cathy O'Neill \cite{Cat} warns us that big data and algorithms, under the impressive but fake mantel of objective truth that mathematics endows them with, can have significant impact on society. She is no luddite, and she knows that algorithm-based decision-making is likely here to stay. But she warns us nonetheless. Similarly I know that math play will go on as long as humans do mathematics. And I really do not want our joy of doing useless math to die from suffocation due to moral analysis paralysis. But maybe we need to think about our mathematics not only in terms of value judgments represented by words such as ``interesting", ``fun", and ``deep", but also in terms of implications and consequences. 

Mathematicians, especially pure mathematicians, have for too long been able to have their cake and eat it too. We have avoided accepting responsibility for the actions of government agencies and yet continued to seek grants from the very same. We have regularly maintained that our work was useless, inconsequential, when it came to helping or hurting people, communities, and other living creatures, and yet we have also argued that one day our work {\it might} have some useful applications so we as mathematicians are useful and should continue to be employed by the public. In short we have claimed to be useful mathematicians doing useless mathematics. 

One of these days we'll need to be honest with everyone, including ourselves.

\end{document}